\documentclass[11pt]{article}

\usepackage{amssymb,latexsym,amsmath}

\usepackage{graphicx}

\begin{document}

\newcommand{\bfi}{\bfseries\itshape}

\makeatletter

\@addtoreset{figure}{section}

\def\thefigure{\thesection.\@arabic\c@figure}

\def\fps@figure{h, t}

\@addtoreset{table}{bsection}

\def\thetable{\thesection.\@arabic\c@table}

\def\fps@table{h, t}

\@addtoreset{equation}{section}

\def\theequation{\thesubsection.\arabic{equation}}

\makeatother

\newtheorem{thm}{Theorem}[section]

\newtheorem{prop}[thm]{Proposition}

\newtheorem{lema}[thm]{Lemma}

\newtheorem{cor}[thm]{Corollary}

\newtheorem{defi}[thm]{Definition}

\newtheorem{rk}[thm]{Remark}

\newtheorem{exempl}{Example}[section]

\newenvironment{exemplu}{\begin{exempl}  \em}{\hfill $\surd$

\end{exempl}}

\newcommand{\comment}[1]{\par\noindent{\raggedright\texttt{#1}

\par\marginpar{\textsc{Comment}}}}

\newcommand{\todo}[1]{\vspace{5 mm}\par \noindent \marginpar{\textsc{ToDo}}\framebox{\begin{minipage}[c]{0.95 \textwidth}

\tt #1 \end{minipage}}\vspace{5 mm}\par}

\newcommand{\ea}{\mbox{{\bf a}}}

\newcommand{\eu}{\mbox{{\bf u}}}

\newcommand{\ueu}{\underline{\eu}}

\newcommand{\ueo}{\overline{u}}

\newcommand{\oeu}{\overline{\eu}}

\newcommand{\ew}{\mbox{{\bf w}}}

\newcommand{\ef}{\mbox{{\bf f}}}

\newcommand{\eF}{\mbox{{\bf F}}}

\newcommand{\eC}{\mbox{{\bf C}}}

\newcommand{\en}{\mbox{{\bf n}}}

\newcommand{\eT}{\mbox{{\bf T}}}

\newcommand{\eL}{\mbox{{\bf L}}}

\newcommand{\eR}{\mbox{{\bf R}}}

\newcommand{\eV}{\mbox{{\bf V}}}

\newcommand{\eU}{\mbox{{\bf U}}}

\newcommand{\ev}{\mbox{{\bf v}}}

\newcommand{\eve}{\mbox{{\bf e}}}

\newcommand{\uev}{\underline{\ev}}

\newcommand{\eY}{\mbox{{\bf Y}}}

\newcommand{\eK}{\mbox{{\bf K}}}

\newcommand{\eP}{\mbox{{\bf P}}}

\newcommand{\eS}{\mbox{{\bf S}}}

\newcommand{\eJ}{\mbox{{\bf J}}}

\newcommand{\eB}{\mbox{{\bf B}}}

\newcommand{\eH}{\mbox{{\bf H}}}

\newcommand{\leb}{\mathcal{ L}^{n}}

\newcommand{\eI}{\mathcal{ I}}

\newcommand{\eE}{\mathcal{ E}}

\newcommand{\hen}{\mathcal{H}^{n-1}}

\newcommand{\eBV}{\mbox{{\bf BV}}}

\newcommand{\eA}{\mbox{{\bf A}}}

\newcommand{\eSBV}{\mbox{{\bf SBV}}}

\newcommand{\eBD}{\mbox{{\bf BD}}}

\newcommand{\eSBD}{\mbox{{\bf SBD}}}

\newcommand{\ecs}{\mbox{{\bf X}}}

\newcommand{\eg}{\mbox{{\bf g}}}

\newcommand{\paromega}{\partial \Omega}

\newcommand{\gau}{\Gamma_{u}}

\newcommand{\gaf}{\Gamma_{f}}

\newcommand{\sig}{{\bf \sigma}}

\newcommand{\gac}{\Gamma_{\mbox{{\bf c}}}}

\newcommand{\deu}{\dot{\eu}}

\newcommand{\dueu}{\underline{\deu}}

\newcommand{\dev}{\dot{\ev}}

\newcommand{\duev}{\underline{\dev}}

\newcommand{\weak}{\stackrel{w}{\approx}}

\newcommand{\mild}{\stackrel{m}{\approx}}

\newcommand{\strong}{\stackrel{s}{\approx}}

\newcommand{\weakdown}{\rightharpoondown}

\newcommand{\opg}{\stackrel{\mathfrak{g}}{\cdot}}

\newcommand{\opunu}{\stackrel{1}{\cdot}}
\newcommand{\opdoi}{\stackrel{2}{\cdot}}

\newcommand{\opn}{\stackrel{\mathfrak{n}}{\cdot}}

\newcommand{\tr}{\ \mbox{tr}}

\newcommand{\Ad}{\ \mbox{Ad}}

\newcommand{\ad}{\ \mbox{ad}}

\renewcommand{\contentsname}{ }

\title{Bipotentials for  non monotone multivalued operators: fundamental results and applications}

\author{Marius Buliga\footnote{"Simion Stoilow" Institute of Mathematics of the 
Romanian Academy, PO BOX 1-764,014700 Bucharest, Romania, e-mail: 
Marius.Buliga@imar.ro } ,
 G\'ery de Saxc\'e\footnote{Laboratoire de 
  M\'ecanique de Lille, UMR CNRS 8107, Universit\'e des Sciences et 
Technologies de Lille,
 Cit\'e Scientifique, F-59655 Villeneuve d'Ascq cedex, France, 
e-mail: gery.desaxce@univ-lille1.fr} , Claude  Vall\'ee\footnote{Laboratoire de 
M\'ecanique des Solides, UMR 6610, UFR SFA-SP2MI, bd M. et P. Curie, 
t\'el\'eport 2, BP 30179, 86962 Futuroscope-Chasseneuil cedex, 
France, e-mail: vallee@lms.univ-poitiers.fr} }

\date{This version: 01.02.2008}

\maketitle

{\bf MSC-class:} 49J53; 49J52; 26B25

\begin{abstract}
This is a survey of recent results about bipotentials representing multivalued operators. 
The notion of bipotential  is  based on an extension of Fenchel's inequality, with  
several interesting applications related to non associated constitutive laws in  
non smooth mechanics, such as  Coulomb frictional contact  or non-associated 
Dr\"ucker-Prager model in plasticity. 
  
Relations betweeen bipotentials and   Fitzpatrick functions are described.  
Selfdual lagrangians,  introduced and studied by Ghoussoub,  can be seen as bipotentials 
representing maximal monotone operators. We show that bipotentials can represent some 
monotone but not maximal operators, as well as non monotone operators. 
 
Further we describe results concerning the construction of a bipotential which represents a given   
non monotone operator, by using convex lagrangian covers or bipotential convex covers.
At the end we prove a new  reconstruction theorem for a bipotential  from a convex 
lagrangian cover, this time using a  convexity notion related to   a minimax theorem of Fan. 
\end{abstract}

\maketitle

\newpage

\tableofcontents

\section{Introduction}

In the generalized standard material theory of Halphen and  Son \cite{Halp JM 75}   
any  constitutive law of a standard material   relates a strain rate variable 
$x \in X$ with a stress-like variable $y \in Y$ by using a dissipation potential $\phi$. 

 $X$ and $Y$ are topological, locally convex, real vector spaces of dual 
variables $x \in X$ and $y \in Y$, with the duality product 
$\langle \cdot , \cdot \rangle : X \times Y \rightarrow \mathbb{R}$. 
 The  dissipation potential $\phi$  is a convex and lower semicontinuous function 
 defined on $X$ and the associated constitutive law is given  by one the following 
 equivalent conditions: 
\begin{enumerate}
 \item[(a)] $\displaystyle y \in \partial \phi(x)$, where $\partial \phi$ is the 
 subdifferential of $\phi$ in the sense of Convex Analysis, 
\item[(b)] $\displaystyle x \in \partial \phi^{*}(y)$, where $\phi^{*}$ is the Fenchel 
dual of $\phi$ with respect to the duality product, 
\item[(c)] $\phi(x) + \phi^{*}(y) = \langle x , y \rangle$ . 
\end{enumerate}

The constitutive laws of standard materials are also called associated laws. From the 
mathematical viewpoint such laws are cyclically monotone operators.  
However, there are many   non-associated constitutive laws, described by a multivalued 
 operator $\displaystyle T:X\rightarrow 2^{Y}$ {\bf which is not cyclically monotone, 
 in some cases not even monotone}.

A possible way to study  non-associated constitutive laws by using 
convex analysis,   proposed first in \cite{saxfeng}, consists in constructing a 
"bipotential" function $b$ of two variables, which physically represents the dissipation.  
See definition \ref{def2} and the section \ref{secbipo} for the introduction into the 
subject  of bipotentials. 

A bipotential  function $b$  is bi-convex, satisfies an inequality generalizing 
Fenchel's one,  namely $ \forall x \in X, y \in Y, \: \ b(x,y) \geq \langle x, y \rangle $, 
and a relation involving partial subdifferentials of $b$ with respect to  variables  
$x$, $y$. 
In the case of associated constitutive laws  the 
bipotential has the expression $b(x,y) = \phi(x)+\phi^{*}(y)$ (section \ref{secbipo}). 
The graph of a bipotential $b$ is simply the set $M(b) \subset X \times Y$ of 
those pairs $(x,y)$ such that $ b(x,y) =  \langle x, y \rangle $.  

A maximal cyclically monotone  operator  $\displaystyle T:X\rightarrow 2^{Y}$ is represented  
by a lower semicontinuous and convex "potential" function $\phi$, by a well known theorem of 
Rockafellar. More general, we say that 
 the bipotential $b$ represents the multivalued operator $T$ 
 if the graph of $T$  equals $M(b)$.

There are already many applications of bipotentials to mechanics. Among them  we cite:  
Coulomb's friction law \cite{sax CRAS 92}, surveyed here in section \ref{secoulomb}, 
non-associated Dr\"ucker-Prager \cite{sax boush KIELCE 93} (section \ref{sedrucker})  
and Cam-Clay models \cite{sax BOSTON 95}  in soil mechanics, cyclic plasticity 
(\cite{sax CRAS 92},\cite{bodo sax EJM 01}) and viscoplasticity \cite{hjiaj bodo CRAS 00} 
of metals with non linear kinematical hardening rule, Lemaitre's damage law \cite{bodo}, 
the coaxial laws (\cite{sax boussh 2},\cite{vall leri CONST 05}), details in sections 
\ref{hillb} and \ref{secoax}.

The notion of a bipotential associated to a multivalued operator is interesting 
also from a mathematical point of view. We show that bipotentials are related 
to  Fitzpatrick functions associated to a maximally monotone operator. Selfdual lagrangians 
introduced and studied by Ghoussoub \cite{ghoussoub1} can be seen as bipotentials 
representing maximal monotone operators (section \ref{fitzp}). In section 
\ref{nonmax}  we describe a result from \cite{bipo3} which implies that some monotone 
non maximal operators can be represented by bipotentials. 

Other  examples of bipotentials come from inequalities. For example, 
 the Cauchy-Bunyakovsky-Schwarz inequality can be recast as: if $X = Y$ is a Hilbert space 
then the function $b: X \times X \rightarrow \mathbb{R}$ defined by  $b(x,y) = \|x\| \|y\|$ 
is a bipotential. More general, some inequalities involving eigenvalues of real 
symmetric matrices can be put in a similar form, thus providing more non trivial 
examples of bipotentials.

In order to better understand the bipotential approach, in  \cite{bipo1}, \cite{bipo3} 
we solved two key problems: (a) when the graph of a given multivalued operator 
can be expressed as the set of critical points of a  bipotentials, and (b) a method 
of construction of a bipotential associated (in the sense of point (a)) to a 
multivalued, typically non monotone, operator. 

At the end of this paper we prove a a new  reconstruction theorem for a bipotential  
from a convex lagrangian cover, this time using a  convexity notion related to   
a minimax theorem of Fan.

\section{Notations and first definitions}

$X$ and $Y$ are topological, locally convex, real vector spaces of dual 
variables $x \in X$ and $y \in Y$, with the duality product 
$\langle \cdot , \cdot \rangle : X \times Y \rightarrow \mathbb{R}$. 
The topologies of the spaces $X, Y$ are  compatible with the duality 
product, that is: any  continuous linear functional on $X$ (resp. $Y$) 
has the form $x \mapsto \langle x,y\rangle$, for some $y \in Y$ (resp. 
$y \mapsto \langle x,y\rangle$, for some  $x \in X$).

We use the notation: $\displaystyle \bar{\mathbb{R}} = \mathbb{R}\cup \left\{ +\infty \right\}$. 

Given a function 
$\displaystyle \phi: X \rightarrow  \bar{\mathbb{R}}$, the domain 
$dom \, \phi$ is the set of points $x \in X$ with $\phi(x) \in \mathbb{R}$.  
The polar of $\phi$, or Fenchel conjugate, 
$\phi^{*}: Y \rightarrow \bar{\mathbb{R}}$ is defined by: 
$\displaystyle \phi^{*}(y) = \sup \left\{ \langle y,x\rangle - 
\phi(x) \mid x \in X \right\}$.  
 
We denote by $\Gamma(X)$ the class of convex and lower semicontinuous 
functions $\displaystyle \phi: X \rightarrow \bar{\mathbb{R}}$. 
The  class of functions $\phi \in \Gamma(X)$  with 
$\displaystyle dom \, \phi \not = \emptyset$ is denoted by 
$\displaystyle \Gamma_{0}(X)$. The class of convex and lower semicontinuous 
$\displaystyle \phi: X \rightarrow \mathbb{R}$ is denoted by 
$\displaystyle \Gamma(X, \mathbb{R})$.

For any convex and closed set $A \subset X$, its  indicator function,  
$\displaystyle \chi_{A}$, is defined by 
$$\chi_{A} (x) = \left\{ \begin{array}{ll}
0 & \mbox{ if } x \in A \\ 
+\infty & \mbox{ otherwise } 
\end{array} \right. $$

\begin{defi} 
The {\bf graph of an operator}   $\displaystyle T: X \rightarrow 2^{Y}$ is the set $G(T)$: 
$$G(T) = \left\{ (x, y) \in X  \times Y \mbox{ : } y \in T(x) \right\}$$
To a graph  $M \subset X \times Y$ we associate the multivalued operators: 
$$\displaystyle  X \ni x  \mapsto m(x) \ = \ \left\{ y \in Y \mid (x,y) \in 
M \right\} \ ,$$
$$\displaystyle  Y \ni y  \mapsto m^{*}(y) \ = \ \left\{ x \in X \mid (x,y) \in M \right\} \ .$$
The {\bf domain of the graph} $M$  is   $\displaystyle dom(M) = \left\{ x \in 
X \mid m(x) \not = \emptyset\right\}$. 
The {\bf image  of the graph}  $M$  is the set $\displaystyle im(M) = \left\{ 
y \in Y \mid m^{*}(y) \not = \emptyset\right\}$. 
\label{def1}
\end{defi}

\begin{defi}
The {\bf subdifferential}   of a function 
$\displaystyle \phi: X \rightarrow \bar{\mathbb{R}}$ in a point 
$x \in dom \, \phi$ is the (possibly empty) set: 
$$\partial \phi(x) = \left\{ u \in Y \mid \forall z \in X  \  
\langle z-x, u \rangle \leq \phi(z) - \phi(x) \right\} \  .$$ 
In a similar way is defined the subdifferential  of a function 
$\psi: Y \rightarrow \bar{\mathbb{R}}$ in a point $y \in dom \, \psi$, 
as the set: 
$$\partial \psi(y) = \left\{ v \in X \mid \forall w \in Y  \  \langle v, w-y \rangle \leq \psi(w) - \psi(y) \right\} \ .$$ 
\end{defi}

With these notations and definitions we have the Fenchel inequality. 

\begin{thm}
Let $\displaystyle \phi: X \rightarrow \bar{\mathbb{R}}$ be a convex lower 
semicontinuous function. Then: 
\begin{enumerate}
\item[(i)] for any $x \in X , y\in Y$ we have $\displaystyle \phi(x) + \phi^{*}(y)  \geq \langle x, y \rangle$; 
\item[(ii)]  for any $(x,y) \in X \times Y$ we have the equivalences: 
$$ y \in \partial \phi(x) \ \Longleftrightarrow \ x \in \partial \phi^{*}(y)  \ \Longleftrightarrow \  \phi(x) + \phi^{*}(y)  = 
\langle x , y \rangle \ . $$
\end{enumerate}
\end{thm}

\begin{defi}
An operator  $\displaystyle T: X \rightarrow 2^{Y}$ is {\bf monotone} if for any $(\displaystyle (x, y), (x', y') \in G(T)$ we have 
$$\langle x - x' , y - y' \rangle \geq 0$$
A graph $M \subset X \times Y$ is monotone if it is the graph of a monotone operator. 
The graph is {\bf maximal monotone} 
(or the associated operator is maximal monotone) if for any monotone graph $M' \subset X\times Y$ such that $M \subset M'$ 
we have $M = M'$. 

An operator  $\displaystyle T: X \rightarrow 2^{Y}$  is cyclically monotone if its graph $G(T)$ is cyclically monotone. A 
graph $M$ is  {\bf cyclically monotone} if for all integer $m>0$ and any finite family of couples $(x_{j},y_{j}) \in M, j=0,1,\ldots,m$, 
\begin{equation}
\displaystyle  \langle x_{0} -  x_{m}, y_{m} \rangle + \sum_{k=1}^{m} \langle x_{k} -  x_{k-1}, y_{k-1} \rangle \leq 0.
\label{Cyclically Monotone inequality}
\end{equation}
A cyclically monotone graph $M$ is  maximal if it does not admit a strict prolongation which is 
cyclically monotone.
\end{defi}

\section{Bipotentials}
\label{secbipo}

\begin{defi} A {\bf bipotential} is a function $b: X \times Y \rightarrow 
\bar{\mathbb{R}}$ with the properties: 
\begin{enumerate}
\item[(a)] $b$ is convex and lower semicontinuous in each argument; 
\item[(b)] for any $x \in X , y\in Y$ we have $\displaystyle b(x,y) \geq 
\langle x, y \rangle$; 
\item[(c)]  for any $(x,y) \in X \times Y$ we have the equivalences: 
\begin{equation}
y \in \partial b(\cdot , y)(x) \ \Longleftrightarrow \ x \in \partial 
b(x, \cdot)(y)  \ \Longleftrightarrow \ b(x,y) = 
\langle x , y \rangle \ .
\label{equiva}
\end{equation}
\end{enumerate}
The {\bf graph of} $b$ is 
\begin{equation}
M(b) \ = \ \left\{ (x,y) \in X \times Y \ \mid \ b(x,y) = \langle x, y 
\rangle \right\} \  .
\label{mb}
\end{equation}
\label{def2}
\end{defi}

An equivalent definition of a bipotential comes out from the following proposition. 

\begin{prop}
 A function $b: X \times Y \rightarrow 
\bar{\mathbb{R}}$ is a bipotential if and only if the following conditions are satisfied: 
\begin{enumerate}
 \item[(A)] $b$ is convex and lower semicontinuous in each argument and  for any $x \in X , y\in Y$ we have $\displaystyle b(x,y) \geq 
\langle x, y \rangle$; 
\item[(B1)] for any $y \in Y$, if  the function $\displaystyle z \in X \mapsto \left( b(z,y) - \langle z, y \rangle \right)$ 
has a minimum, then the minimum equals $0$; 
\item[(B2)] for any $x \in X$, if  the function $\displaystyle p \in Y \mapsto \left( b(x,p) - \langle x, p \rangle \right)$ 
has a minimum, then the minimum equals $0$.
\end{enumerate}
\label{punus}
\end{prop}

\paragraph{Proof.} 
Condition (A) is the same as conditions (a),(b) from definition \ref{def2}. All we have to prove is: if the function $b$ satisfies 
condition (A) then condition  (c) from definition \ref{def2} is equivalent with (B1) and (B2). 

Assume (A) and take $x \in X$, $y \in Y$ such that  $x \in \partial 
b(x, \cdot)(y)$. This is equivalent with: $x$ is a minimizer of the function 
$$\displaystyle z \in X \mapsto \left( b(z,y) - \langle z, y \rangle \right)$$ 
Therefore the statement  $x \in \partial 
b(x, \cdot)(y)  \ \Longleftrightarrow \ b(x,y) = 
\langle x , y \rangle$ is equivalent with (B1). In the same way we prove that $y \in \partial b(\cdot , y)(x)  \ \Longleftrightarrow \ b(x,y) = 
\langle x , y \rangle$ is equivalent with (B2).
\quad $\square$

This simple proposition justifies the introduction of strong bipotentials, which are particular cases of bipotentials. Conditions 
(B1S) and (B2S) appeared first time as relations (51), (52) \cite{laborde}, 

\begin{defi}
 A function $b: X \times Y \rightarrow 
\bar{\mathbb{R}}$ is a {\bf strong bipotential} if it satisfies the conditions: 
\begin{enumerate}
\item[(a)] $b$ is convex and lower semicontinuous in each argument; 
\item[(B1S)] for any  $y \in Y$ $\inf \left\{  b(z,y) - \langle z, y \rangle \mbox{ : } z \in X \right\} \in  \left\{ 0, + \infty \right\}$; 
\item[(B2S)] for any  $x \in X$ $\inf \left\{  b(x,p) - \langle x, p \rangle \mbox{ : } p \in Y \right\} \in  \left\{ 0, + \infty \right\}$. 
\end{enumerate}
\label{defbipos}
\end{defi}

\section{Operators representable by a bipotential}

\begin{defi}  The non empty set $M \subset X \times Y$ is a  {\bf BB-graph}  
(bi-convex, bi-closed) if for all $x \in \ dom(M)$ and for all $y \in \ im(M)$ 
the sets $\displaystyle m(x)$ and $\displaystyle m^{*}(y)$ are convex and closed.
\label{dh1}
\end{defi}

The following theorem (theorem 3.2 \cite{bipo1}) gives a necessary and   sufficient condition for the 
existence of a bipotential associated to a constitutive law $M$. 

\begin{thm}
 Given  a non empty set $M \subset X \times Y$, there is a bipotential 
$b$ such that $M=M(b)$ if and only if $M$ is a BB-graph. 
 \label{thm1}
 \end{thm}

The bipotential mentioned in the previous theorem is denoted by $\displaystyle b_{M}$ and it has the expression: 
\begin{equation}
b_{M} (x,y) \ = \ \langle x , y \rangle + \chi_{M}(x,y) 
\label{binfty}
\end{equation}

In the case of a maximal cyclically monotone 
graph $M$, by Rockafellar theorem (\cite{rocka} Theorem 24.8.) there is an unique separable
bipotential associated to $M$ ( see section \ref{sepbipo} for the
definition of separable bipotentials). With the bipotential given by (\ref{binfty}) we have 
two different bipotentials representing the same graph. Therefore, in the larger class 
made of all bipotentials, in general there is no unicity of the bipotential representing a 
given BB-graph.

For any BB-graph $M$, a bipotential $b$ is admissible if $M \subset M(b)$. Then we obviously have 
$\displaystyle b(x,y) \leq b_{M}(x,y)$ for any $(x,y) \in X \times Y$. In this sense $\displaystyle 
b_{M}$ is the greatest admissible bipotential for the BB-graph $M$.

\section{Fitzpatrick functions; selfdual lagrangians} 
\label{fitzp}

Let $X$ be a reflexive Banach space 
and $\displaystyle X^{*}$ its topological dual. The duality product between $X$ and $\displaystyle X^{*}$  is the function 
$\pi: X \times X^{*} \rightarrow \mathbb{R}$, defined by $\displaystyle \pi(x,x^{*}) = \langle x, x^{*} \rangle = x^{*}(x)$. 

The space $\displaystyle X \times X^{*}$ is in duality with itself  by the duality product 
$$ \langle (x,x^{*}) , (y, y^{*}) \rangle = \langle x, y^{*} \rangle +  \langle y, x^{*} \rangle $$

Fitzpatrick functions have been introduced in \cite{fitzpatrick1}. More on Fitzpatrick functions can be found in 
\cite{borwein1} and in the book \cite{borwein2}. 
\begin{defi}
 The {\bf Fitzpatrick function}  associated to a graph $M \subset X \times X^{*}$  is the  function 
$\displaystyle f_{M}: X \times X^{*} \rightarrow \mathbb{R} \cup \left\{ + \infty \right\}$ given by the Fenchel dual of 
$\displaystyle b_{M}$. Equivalently, $\displaystyle f_{M}$ is given by: 
$$f_{M}(x,x^{*}) = \sup \left\{ \langle a, x^{*}\rangle + \langle x, a^{*}\rangle -  \langle a, a^{*}\rangle \mbox{ : } (a, a^{*}) \in M \right\}$$
\end{defi}

\begin{prop}
 (Properties of the Fitzpatrick function) Let $M \subset X \times X^{*}$ be a  graph. Then the associated Fitzpatrick function 
$\displaystyle f_{M}$ has the properties: 
\begin{enumerate}
\item[(a)] $\displaystyle f_{M}$ is convex and lower semicontinuous, 
 \item[(b)] the graph $M$ is maximal monotone if and only if: 
\begin{enumerate}
\item[(b1)]  for any $(x,x^{*}) \in X \times X^{*}$ we have $\displaystyle f_{M}(x,x^{*}) \geq \langle x, x^{*} \rangle$ and  
\item[(b2)] we have equality  $\displaystyle f_{M}(x,x^{*}) = \langle x, x^{*} \rangle$ if and only if $(x, x^{*}) \in M$. 
\end{enumerate}
\item[(c)] if $M$ is a BB-graph then the function 
$$g_{M}(x,x^{*}) \ = \ \left\{ \begin{array}{lcl}
                                f_{M}(x,x^{*}) & , & (x, x^{*}) \in \, dom \, M \, \times \, im \, M \\
                                +\infty  &  & \mbox{ otherwise}  
                               \end{array} \right.$$ 
 is a strong bipotential.
\end{enumerate}
\label{propfitz}
\end{prop}

\paragraph{Proof.}
 By construction any Fitzpatrick function is convex and lower semicontinuous. For proving (b) it is enough to use the following characterization of 
maximal monotone graphs: $M$ is maximal monotone if and only if 
\begin{enumerate}
 \item[(i1)] for any $(x,x^{*}) \in X \times X^{*}$ we have 
$$\displaystyle \inf \left\{ \langle x - a , x^{*} - a^{*} \rangle \mbox{ : } (a, a^{*}) \in M \right\} \leq 0$$
\item[(i2)] $(x,x^{*}) \in M$ if and only if 
$$\displaystyle \inf \left\{ \langle x - a , x^{*} - a^{*} \rangle \mbox{ : } (a, a^{*}) \in M \right\} = 0$$
\end{enumerate}
The Fitzpatrick function associated to the graph $M$ can be written as: 
$$f_{M}(x, x^{*}) \ = \ \langle x, x^{*}\rangle \, - \, \inf \left\{ \langle x - a , x^{*} - a^{*} \rangle \mbox{ : } (a, a^{*}) \in M \right\}$$
Therefore (b1), (b2) follow from (i1), (i2) respectively. 

In order to prove (c) we have to check (B1S), (B2S) definition \ref{defbipos}.  Let $x \in \, dom \, M$. Then 
$$ \inf \left\{  g_{M}(x,x^{*}) - \langle x, x^{*} \rangle \mbox{ : } x^{*} \in \, im \, M \right\}  \, = \,  $$
$$ = \, - \inf \left\{ \langle x - a , x^{*} - a^{*} \rangle \mbox{ : } (a, a^{*}) \in M \, , \, x^{*} \in \, im \, M \right\} \ = \ 0 $$
The proof of (B2S) is similar. 
\quad $\square$

Further we describe selfdual lagrangians, with the notations from  Ghoussoub \cite{ghoussoub1}. See 
the mentioned paper and the references therein, especially \cite{ghoussoub2} \cite{ghoussoub3}
 for more informations on the variational theory associated to selfdual lagrangians. 
 Here we point out that selfdual lagrangians are 
particular cases of bipotentials.

\begin{defi}
Let $X$ be a reflexive Banach space.  To any function $L \in \Gamma_{0}(X \times X^{*})$ we associate the following operators: 
\begin{enumerate}
 \item[(i)] $\displaystyle \delta L : X \rightarrow 2^{X^{*}}$ defined by: 
$$\delta L (x) \ = \ \left\{ x^{*} \in X^{*} \mbox{ : } (x, x^{*}) \in \partial L (x, x^{*}) \right\} $$ 
(here $\partial L (x, x^{*})$ is the subdifferential of $L$), 
\item[(ii)] $\displaystyle \bar{\partial} L : X \rightarrow 2^{X^{*}}$ defined by: 
$$\bar{\partial} L (x)  \ = \ \left\{ x^{*} \in X^{*} \mbox{ : } L(x, x^{*} \ = \ \langle x , x^{*} \rangle \right\}$$
\end{enumerate}

\end{defi}

\begin{defi}
 A function $L \in \Gamma_{0}(X \times X^{*})$ is a {\bf selfdual lagrangian} if $L = L^{*}$. 
\end{defi}

The following is  a slight reformulation of lemma 2.1 and proposition 2.1 \cite{ghoussoub1} 

\begin{prop}
 If $L$ is a selfdual lagrangian such that for some 
$\displaystyle x_{0} \in X$ the function $\displaystyle L(x_{0}, \cdot )$ is bounded on the balls of $X^{*}$, 
then $L$   is a strong  bipotential and we have $M(L) = G(\bar{\partial} L) = G(\delta L)$.  Moreover, in this case for any $x^{*} \in X^{*}$ there exists $\bar{x} \in X$ such that $x^{*} \in \delta L(\bar{x})$ and 
$$L(\bar{x}, x^{*}) \ - \ \langle \bar{x} , x^{*} \rangle \ = \ \inf \left\{ L(x, x^{*}) - \langle x , x^{*} \rangle \mbox{ : } x \in X \right\} \ = \ 0$$
\end{prop}

\paragraph{Proof.}
 As mentioned before, from lemma 2.1 \cite{ghoussoub1} we get that for any selfdual lagrangian $L$ we have $\bar{\partial} L = \delta L$. By definition 
any selfdual lagrangian is convex, lower semicontinuous and for any $(x, x^{*}) \in X \times X^{*}$ we have $L(x,x^{*}) \geq \langle x, x^{*} \rangle$, 
as a consequence of the Fenchel inequality in $X \times X^{*}$. The fact that $L$ is a strong bipotential (conditions (BS1), (BS2) 
definition \ref{defbipos}), as well as the final part of the conclusion are straightforward reformulations of the conclusion of 
 Proposition 2.1 \cite{ghoussoub1}. 
\quad $\square$

The following proposition is an application of lemma 3.1 and proposition 3.1 \cite{ghoussoub1}.

\begin{prop}
 Let $\displaystyle M \subset X \times X^{*}$ be a maximal monotone graph. Then there exists a selfdual 
 lagrangian $L_{M}$ such that $\displaystyle G(\bar{\delta} L_{M}) = M$. 
\end{prop}

\paragraph{Proof.}
 Let $\displaystyle f_{M}$ be the Fitzpatrick function of $M$. Then, according to lemma 3.1 \cite{ghoussoub1} the Fitzpatrick function 
$\displaystyle f_{M}$ satisfies the hypothesis of proposition 3.1 \cite{ghoussoub1}. Therefore the selfdual lagrangian defined by
$$L_{M}(x,x^{*}) \ = \ \inf \left\{ \frac{1}{2} f_{M}(x_{1}, x_{1}^{*}) + \frac{1}{2} f_{M}^{*}(x_{2}, x_{2}^{*}) + \frac{1}{8} \| x_{1} - x_{2} \|^{2} + 
\right. $$
$$\left. +  \frac{1}{8} \| x_{1}^{*} - x_{2}^{*} \|^{2} \mbox{ : } (x, x^{*}) = \frac{1}{2} (x_{1}, x_{1}^{*}) + \frac{1}{2} (x_{2}, x_{2}^{*}) \right\}$$ 
mentioned in the proof of the  proposition 3.1 \cite{ghoussoub1} as ``the proximal average'' between $\displaystyle f_{M}$ and $\displaystyle f_{M}^{*}$ 
is the one needed. 
\quad $\square$

\section{Separable bipotentials}
\label{sepbipo}
 If 
$\phi: X \rightarrow \mathbb{R}$ is a convex, lower semicontinuous  potential, 
consider the multivalued operator $\displaystyle \partial \phi$ (the 
subdifferential of $\phi$).  The graph of this operator  is the set 
\begin{equation}
M(\phi) \ = \ \left\{ (x,y) \in X \times Y \ \mid \ \phi(x)+\phi^{*}(y) = \langle x, y \rangle \right\} \  .
\label{mphi}
\end{equation}
$M(\phi)$ is  maximal  cyclically 
monotone \cite{rocka} Theorem 24.8. Conversely, if 
$M$ is closed and maximally cyclically monotone then there is a convex, 
lower semicontinuous $\phi$ such that $M=M(\phi)$.  

To  the function  $\phi$  we  associate the 
{\bf separable bipotential} $$\displaystyle b(x,y)  = \phi(x) + \phi^{*}(y) .$$ 
Indeed, the Fenchel inequality can be reformulated by saying that the function 
$b$,  previously defined, is a bipotential. More precisely, the point (b) 
(resp. (c))  in the definition of a bipotential corresponds to (i) 
(resp. (ii)) from Fenchel inequality.

The bipotential $b$ and the function  
$\phi$ have the same graph, that is  $\displaystyle M(b) = M(\phi)$.

\section{Bipotentials for monotone, non maximal graphs}
\label{nonmax}

The following two results are from the paper \cite{bipo3} (theorem 3.1 and 
corollary 3.3).

In the theorem below it is shown that intersections of two maximal monotone 
graphs are sometimes representable by a bipotential. Therefore there exist bipotentials 
$b$ with $M(b)$ monotone, but not maximal. 

\begin{thm}
Let  $b_{1}$ and $b_{2}$ be separable bipotentials associated respectively to 
the convex and lower semicontinuous functions 
$\displaystyle \phi_{1} , \phi_{2}: X
\rightarrow \mathbb{R} \cup \left\{ + \infty \right\}$, that is 
$$b_{i} (x,y) \ = \ \phi_{i}(x) + \phi_{i}^{*}(y)$$ 
for any $i = 1,2$ and $(x,y) \in X \times Y$. Consider the following 
assertions:
\begin{enumerate}
\item[(i)] 
$b= max (b_{1},b_{2})$ is a strong 
bipotential and $M(b)=M(b_{1})\cap M(b_{2})$.
\item[(ii')] For any $\displaystyle y \in \, dom \,  \phi_{1}^{*} \, \cap \, dom \,  \phi_{2}^{*}$ and for any $\lambda \in 
[0,1]$ we have 
\begin{equation}
\left( \lambda \, \phi_{1} \, + \, (1-\lambda) \, \phi_{2} \right)^{*}(y) \ = 
\ \lambda \, \phi_{1}^{*}(y) \, + \, (1-\lambda) \, \phi_{2}^{*}(y) 
\label{iiprim}
\end{equation}
\item[(ii")] For any $\displaystyle x \in \, dom \,  \phi_{1} \, \cap \, dom \,  \phi_{2}$ and for any $\lambda \in 
[0,1]$ we have 
\begin{equation}
\left( \lambda \, \phi_{1}^{*} \, + \, (1-\lambda) \, \phi_{2}^{*} 
\right)^{*}(x) \ = 
\ \lambda \, \phi_{1}(x) \, + \, (1-\lambda) \, \phi_{2}(x)
\label{iisec}
\end{equation}
\end{enumerate} 
Then the point (i) is equivalent with the conjunction of (ii'), (ii"), (for
short: 
  (i) $\Longleftrightarrow$ ( (ii') AND (ii") ) ). 
\label{thmsup}
\end{thm}

In the proof of this theorem we make use of a minimax result by Sion \cite{sion}. Notice that 
in section \ref{newconst},  we use another minimax result of Fan \cite{kyfan} in the proof of 
theorem \ref{thm4}.

The conditions (ii'), (ii") from Theorem \ref{thmsup} imply relations which 
can be expressed with the help of inf convolutions. Consider $\displaystyle \phi_{1},
\phi_{2} \in \Gamma_{0}(X)$. For any $\lambda \in (0,1)$ we introduce two
functions defined on $X$ by: 
$$f_{1, \lambda} (x) \, = \, \lambda \, \phi_{1}(\frac{1}{\lambda} x) \quad , 
\quad f_{2, \lambda} (x) \, = \, (1-\lambda) \, \phi_{2}(\frac{1}{1-\lambda} x)$$

\begin{prop}
Let $\displaystyle \phi_{1},
\phi_{2} \in \Gamma_{0}(X)$ such that 
$$b(x,y) \, = \, \max \left(\phi_{1}(x) + \phi_{1}^{*}(y), \phi_{2}(x) +
\phi_{2}^{*}(y) \right\}$$ 
is a strong bipotential. Then, with the previous notations, for any  
$\displaystyle x \in \, dom \,  \phi_{1} \, \cap \, dom \,  \phi_{2}$ and for 
any $\lambda \in (0,1)$ the subdifferential of the inf-convolution 
$\displaystyle f_{1,\lambda} \square f_{2,\lambda}$ has the expression: 
$$\partial \left( f_{1,\lambda} \square f_{2,\lambda} \right) (x) \, = \, 
\partial \phi_{1} (x) \, \cap \, \partial \phi_{2} (x)$$
\label{pconvo}
\end{prop}

\section{Bipotentials and inequalities}   

A source of interesting bipotentials is provided by inequalities. Here we discuss about 
the Cauchy-Bunyakovsky-Schwarz inequality and about an inequality of Fan concerning 
eigenvalues of symmetric matrices. 

\subsection{Cauchy bipotential} 
Let $X=Y$ be a Hilbert space and let the duality 
product be  equal 
to the scalar product. Then we define the {\bf Cauchy bipotential} by the formula 
$$\displaystyle b(x,y) = \| x\| \  \| y\| .$$ 
Let us check the Definition (\ref{def2}) The point (a) is obviously satisfied. 
The point (b) is true by the Cauchy-Bunyakovsky-Schwarz inequality.  
We have equality in the Cauchy-Bunyakovsky-Schwarz inequality 
$b(x,y) = \langle x,y \rangle$ if and only if there is $\lambda > 0$ such 
that $y = \lambda x$ or one of $x$ and $y$ vanishes.  This is exactly the statement from the  point (c), for 
the function $b$ under study.

The graph of the Cauchy bipotential is not monotone.

\subsection{Hill bipotential}
\label{hillb}

 Let $\displaystyle S(n)$ be the space of $n \times n$  real symmetric matrices. There is a bipotential expressing that 
two matrices $X$ and $Y$ have a simultaneous ordered spectral decomposition, \cite{vall leri CONST 05}. In Mechanics, a constitutive law between two 
tensors implying that they admit the same eigenvectors is said to be coaxial \cite{sax boussh 2}, \cite{vall leri CONST 05}.

The space $S(n)$ is endowed with the scalar product 
$$\langle X, Y \rangle = tr \left(X Y \right)$$ 
 Consider the function 
$b: S(n) \times S(n) \rightarrow \mathbb{R}$ with the expression 
$$b(X,Y) \ = \ \lambda_{1}(X) \lambda_{1}(Y) + ... + \lambda_{n}(X) \lambda_{n}(Y)$$ 
where for any $X \in S(n)$  $\displaystyle \lambda_{1}(X) \geq \lambda_{2}(X) \geq ... \geq \lambda_{n}(X)$ are the eigenvalues of $X$ ordered from the largest to the smallest. 
With the help of this function we can write one of Fan's inequalities \cite{fan1} as: for any $X, Y \in S(n)$ we have 
$$b(X,Y) \geq \langle X, Y \rangle$$
with equality if and only if $X$ and $Y$ have the same eigenvectors with preservation of the order of the eigenvalues. 

In \cite{vall leri CONST 05} is proved that for $n=3$ the function $b$ is a bipotential, called {\bf
Hill bipotential} due to applications dealt with by Hill in mechanics. A similar  proof, involving majorization,  can be done for the case of a general $n$. The graph of the Hill bipotential is not monotone.

\section{Bipotentials in non smooth mechanics}

A simple example of a monotone operator in non smooth mechanics   is provided by the following model of plasticity of metals. $X=Y$ is the space of $n \times n$ real symmetric traceless matrices with the pairing $\langle x, y \rangle = \, tr \, xy$  and the associated norm $ \| x \| = \mid \left\langle x,y \right\rangle  \mid^{\frac{1}{2}} $. Let $c$ be a non negative constant. The plasticity operator is defined by 
$$\displaystyle T_{p}  \left(  0 \right) 
   = K = \left\lbrace y \in Y \mid \ \ \| y \| \ \leq c \ \right\rbrace $$
otherwise
$$ \displaystyle T_{p} \left(  x \right)  = \frac{x}{ \|  x \| } $$
In plasticity, $x$ is the plastic strain rate tensor, $y$ is the deviatoric stress tensor and $c$ is the yield stress. The closed convex set $K$ 
is the {\bf plastic domain} and the irreversible or plastic deformations varies when $x\neq0$. The plastic model is not limited to the metals but can be used also for soil materials. 

A more involved  example is the {\bf associated Dr\"ucker-Prager model} where the variable $x$ is the
plastic strain rate tensor and $y$ is the stress tensor (both seen as elements of $S(n)$). 
With usual notations, the tensors $x$ and $y$ are split into their deviatoric and spheric parts:
$$\displaystyle x = x_{d} + \frac{1}{3} x_{h} I,
 \quad  y = y_{d} + y_{h} I   $$
where $I$ is the identity operator, $x_{h} = tr \left( x \right) $ and $y_{h} = \frac{1}{3} tr \left( y \right) $. As the decomposition is unique, we can write in short $x = \left( x_{d}, x_{h} \right) $, $y = \left( y_{d}, y_{h} \right) $ and the duality pairing is $ \left\langle  x,y \right\rangle = tr \left( x_{d} y_{d} \right)  + x_{h} y_{h} $. 

The convex cone parameterized by the
friction angle $\varphi \in  \left( 0,\frac{\pi}{2}\right) $, the cohesion stress $ c > 0 $, and
$r=3\sqrt{2}/\sqrt{9+12\tan^2\varphi}$ (\cite{boush chaa IJP 01}), with the vertex at $\displaystyle v = \frac{c}{\tan\varphi}I$, given by 
$$K=\{y \in Y \ |\ ||y_{d}||\leq r(c-y_{h} \tan\varphi)\}$$ 
is called the plastic domain.

The multivalued operator corresponding to the associated  Dr\"ucker-Prager model is defined by: 
$\displaystyle T_{DP}  \left(  0 \right) = K$, 
if $\displaystyle x_{h} = r|| x_{d} ||\tan\varphi$ then 
$$ \displaystyle T_{DP} \left(  x \right)  = \left\lbrace  y \in Y 
     \mid \exists \eta \geq 0, \ y = v + \eta \left( \frac{x_{d}}{ \|  x_{d} \| } 
     - \frac{1}{r c} v \right) \right\rbrace \quad ,  $$
 if $\displaystyle  x_{h} > r|| x_{d} ||\tan\varphi$ then 
$\displaystyle T_{DP}  \left(  x \right) = \left\lbrace v \right\rbrace $,  
otherwise $ \displaystyle T_{DP} \left(  x \right)  = \emptyset $.

\subsection{Dr\"ucker-Prager non associated plasticity}
\label{sedrucker}

The {\bf non associated Dr\"ucker-Prager model} is  characterized,  as previously, by the friction angle
$\varphi \in  \left( 0,\frac{\pi}{2}\right) $, the cohesion stress $ c > 0 $, and
$r=3\sqrt{2}/\sqrt{9+12\tan^2\varphi}$ but also by a new parameter, the dilatancy angle  
$\theta \in [0,\varphi) $. Once again, $X = Y = S(n) $ and we use the splitting into deviatoric and spheric parts. The associated multivalued operator is $\displaystyle T_{na}$ defined by: $\displaystyle T_{na}  \left(  0 \right) = K$, 
if $ x_{h} = r|| x_{d} ||\tan\theta$ then 
$$ \displaystyle T_{na} \left(  x \right)  = \left\lbrace  y \in Y 
     \mid \exists \eta \geq 0, \ y = v + \eta \left( \frac{x_{d}}{ \|  x_{d} \| } 
     - \frac{1}{r c} v \right) \right\rbrace \quad , $$
if $ x_{h} > r|| x_{d} ||\tan\theta$ then $\displaystyle T_{na}  \left(  x \right) = \left\lbrace v \right\rbrace  $ 
otherwise $ \displaystyle T_{na}  \left(  x \right)  = \emptyset $. 

If we put $\theta = \varphi$ then  we recover the operator $T_{DP}$ of the associated 
Dr\"ucker-Prager model defined above.

The non-associated Dr\"ucker-Prager law $\displaystyle y \in T_{na}(x)$ is equivalent with the following differential inclusion: 
$$x \ + \ \frac{1}{3} \left( x_{h} \  + \ r \, \| x \| \, (\tan \phi \ - \ \tan \theta)\right) \ \in \ \partial \chi_{K}(y)$$
According to \cite{boush chaa IJP 03} \cite{sax boussh 2} \cite{hjiaj fort IJES 03}, this inclusion can be written as $\displaystyle b_{p}(x,y) = \langle x, y \rangle$, for the bipotential 
$$b_{p}(x,y) \ = \ \chi_{K}(y) \ + \ \chi_{K_{p}}(x) \ + \ \frac{c}{\tan \phi} x_{h} \ + $$ 
$$+ \ r \, \| x \| \, (\tan \phi \ - \ \tan \theta) \left( 
y_{h} \ - \ \frac{c}{\tan \phi}\right)$$
The last term in this expression is a coupling term which gives the implicit character to the constitutive law.

\subsection{Coulomb's dry friction law}
\label{secoulomb}

Another interesting operator comes in relation with the  unilateral contact with dry
friction or {\bf Coulomb's friction model}. T
Consider two bodies $\Omega_{1}$ and $\Omega_{2}$ which are in contact at a point $M$, with $\mathbf{n}$ the unit vector normal to the common tangent plane and directed towards $\Omega_{1}$. The space $X = \mathbb{R}^{3}$ is the one of relative velocities between points of contact of two bodies, and the space $Y$, identified also to $\mathbb{R}^{3}$, is the one of the contact reaction stresses. The duality product is the usual scalar  product. We put
$$(x_{n},x_{t})\in X = \mathbb{R} \times \mathbb{R}^{2}, \quad 
    (y_{n},y_{t})\in Y = \mathbb{R} \times \mathbb{R}^{2}\ ,$$
where $x_{n}$ is the gap velocity, $x_{t}$ is the sliding velocity, $y_{n}$ is the contact pressure and $y_{t}$ is the friction stress. The friction coefficient is $\mu > 0$. 

The graph of the law of unilateral 
contact with Coulomb's dry friction is  the union of three sets, 
respectively corresponding to the 'body separation', the 'sticking' and the 
'sliding'. 
\begin{equation}
M = \left\{ (x,0)\in X \times Y \   \mid \ x_{n} < 0 \right\} 
     \cup 
\label{Coulomb friction contact law}
\end{equation}
$$ \cup \ \left\{ (0,y)\in X \times Y \   \mid \  
                       \parallel y_{t} \parallel \leq \mu y_{n}  \right\} 
		       \cup  $$ 
$$\cup \ \left\{ (x,y) \in X \times Y \   \mid \  x_{n} = 0, \ x_{t} \neq 0, \ 
                       y_{t} = \mu y_{n} \frac{x_{t}}{\parallel x_{t}
		       \parallel}   \right\} $$
It is well known that this graph is not monotone, then not cyclically 
monotone.

This law can be written as the following differential inclusion (\cite{saxfeng} \cite{sax CRAS 92} \cite{sax feng IJMCM 98} \cite{sax boussh 2}): 
$$\left( x_{n} \ - \ \mu \, \|x_{t} \|\right) \mathbf{n} \ + \ x_{t} \ \in \ \partial \chi_{K_{\mu}}(y)$$

Let us consider the  conjugate cone of the Coulomb's cone: 
$$ K_{\mu}^* =  \left\{ (x_{n},x_{t})\in X \   \mid \  
                      \mu \parallel x_{t} \parallel  + x_{n} \leq 0  \right\}  .$$
We shall use also a second pair of conjugate cones: 
$$ K_{0} =  \left\{ (y_{n},0)\in Y \   \mid \  y_{n} \geq 0  \right\}  , 
\quad 
     K_{0}^* =  \left\{ (x_{n},x_{t})\in X \   \mid \   x_{n} \leq 0  \right\}  .$$

The graph of the law of unilateral 
contact with Coulomb's dry friction is the graph of the following bipotential,  previously given in \cite{saxfeng}:
$$b (x,y) \ =  \mu y_{n} \parallel x_{t} \parallel + \chi_{K_{\mu}} (y) +  \chi_{ K^{*}_{0} } (x) \ .$$
 
\subsection{Coaxial laws}
\label{secoax}

Consider the non monotone operator defined by 
$ \displaystyle T_{iso} \left(  0 \right)  =  \mathbb{R}^{n}$,  
otherwise 
$$ \displaystyle T_{iso} \left(  x \right)  = 
     \left\lbrace y \in Y \mid \exists \lambda > 0, \ \ y = \lambda x \right\rbrace  $$ 
An operator $\displaystyle S: \mathbb{R}^{n} \rightarrow 2^{\mathbb{R}^{n}}$ is strongly 
isotropic if its graph is contained in the graph of $T_{iso}$. The graph of the non monotone 
operator $\displaystyle T_{iso}$ is the graph of the Cauchy bipotential. 

The eigenvalues of any matrix 
$x \in S(n)$ can be conventionally ordered from the largest to the smallest:
$\lambda_1(x)\geq \lambda_2(x)\geq \ldots \geq \lambda_n(x)$. 
Two elements   $x, y \in S(n)$ are said coaxial if they have the same eigenvectors with 
preservation of the order of the eigenvalues: the largest eigenvalue of $x$ is associated to the 
eigenvector of $y$ corresponding to the largest eigenvalue of $y$, and so on. 
An operator $ T: S(n) \rightarrow 2^{S(n)}$ is coaxial (\cite{dangsax},\cite{vall leri CONST 05})
 if for any $y \in T (x)$, $y$ and $x$ are coaxial. 
  The graph of any coaxial  operator 
is  contained in the non monotone operator defined by $ \displaystyle T_{H} \left(  0 \right)  = 
 S(n) $, 
otherwise
$$ \displaystyle T_{H} \left(  x \right)  = 
     \left\lbrace y \in Y \mid x \ \mathrm{and} \ y \ \mathrm{are} \ \mathrm{coaxial} \right\rbrace  $$
We have seen in section \ref{hillb} that the graph of a coaxial operator is included in the graph 
of the Hill bipotential.

\section{Numerical methods based on the bipotential framework}

In applications bipotentials are interesting because of the associated 
implicit normality rules. The properties of bipotentials allow to discretize 
an evolution problem into a series of minimization problems concerning the 
associated {\bf bifunctional}. 

For example let us consider a body with reference configuration $\Omega$, with 
the boundary decomposition $\partial \Omega = \Gamma_{1} \cup \Gamma_{2} \cup \Gamma_{3}$. 
At a fixed moment $t$ we have imposed velocities $\displaystyle \dot{u}_{t}$ on the part 
$\Gamma_{1}$ of its boundary, imposed forces $\displaystyle f_{t}$ on $\Gamma_{2}$ and 
on $\Gamma_{3}$ the body is in unilateral contact with friction with a rigid foundation. The unit 
normal, pointing outwards, of the boundary $\partial \Omega$ is denoted by $\mathbf{n}$ and 
for any stress field $\sigma$ in $\Omega$ we denote  $\sigma_{n} = \sigma \, \mathbf{n}$.

Suppose that the body is made by a plastic material described by a bipotential $\displaystyle 
b_{p}$ dependent on the strain rate $\varepsilon(\dot{u})$ and the stress $\sigma$. The contact with 
friction is described by the bipotential $\displaystyle b_{c}$. associated with the Coulomb law.

For any pair $(v, \tau)$ of kinematically admissible velocity field $v$ and 
statically admissible stress field $\tau$ we introduce the bifunctional
$$B(v, \tau)   \ = \ \int_{\Omega} b_{p}(\varepsilon(v), \tau) \mbox{ d}x \ + \ 
\int_{\Gamma_{3}} b_{c}(-v, \tau_{n}) \mbox{ d}s \ -$$
$$- \ \int_{\Gamma_{1}} \tau_{n} \cdot \dot{u}_{t} \mbox{ d}s \ - \ \int_{\Gamma_{2}} f_{t} \cdot 
v \mbox{ d}s$$
Then, by using integration by parts and the properties of the bipotentials which are involved, 
one can show that 
$$B(v, \tau) \geq 0$$
and that $B(u, \sigma) = 0$ if $(u,\sigma)$ is the pair formed by the velocity field and associated
stress field of the body at moment $t$. 

Therefore we may try to numerically minimize the bifunctional in order to find the solution of 
the (quasistatic) evolution problem. Alternatively, the bifunctional can be adapted 
in order to use a predictor/corrector scheme. 

As applications we may cite the bound theorems of the limit 
analysis (\cite{sax boussh IJMS 98}, \cite{boush chaa IJMS 02}) and the plastic 
shakedown theory (\cite{sax trit AACHEN 00},  \cite{dangsax}, 
\cite{boush chaa IJP 03}, \cite{bouby sax IJSS 06}) which  can be reformulated 
 by means of weak normality rules. The bipotential method suggests new 
algorithms, fast and robust, as well as variational error estimators assessing 
the accurateness of the finite element mesh (\cite{hjiaj sax BUDA 96}, 
\cite{hjiaj sax PARIS 96}, \cite{sax hjiaj GIENS 97},  
\cite{sax feng IJMCM 98}, \cite{boush chaa IJP 01}, \cite{hjiaj fort IJES 03}, 
\cite{hjiaj feng IJNME 04}). Such algorithms have been proposed and used in 
applications to the contact mechanics 
\cite{feng hjiaj CM 06}, the dynamics of granular materials 
(\cite{fort sax CRAS 99},  \cite{fort hjiaj CG 02}, 
\cite{fort mill IJNME 04}\cite{sax fort MOM 04}), the cyclic 
plasticity of metals \cite{sax hjiaj GIENS 97} and the plasticity of 
soils (\cite{bersaxce}, \cite{hjiaj fort IJES 03}).

A very challenging subject seems to be the extension of the mathematical results of 
 Ghoussoub \cite{ghoussoub2} \cite{ghoussoub3} for the particular case of selfdual lagrangians, to 
 general bipotentials, or at least to bipotentials constructed from bipotential convex covers.

\section{Construction of bipotentials}
\label{seconstruct}

Let $Bp(X,Y)$ be the set of all bipotentials $b: X \times Y \rightarrow
\mathbb{R} \cup \left\{ + \infty \right\}$. We shall need the following notion of implicitly convex functions.

\begin{defi}
Let $\Lambda$ be an arbitrary non empty set and $V$ a real vector space. The 
function $f:\Lambda\times V \rightarrow \bar{\mathbb{R}}$ is 
{\bf implicitly  convex} if for any two elements 
$\displaystyle (\lambda_{1}, z_{1}) , 
(\lambda_{2},  z_{2}) \in \Lambda \times V$ and for any two numbers 
$\alpha, \beta \in [0,1]$ with $\alpha + \beta = 1$ there exists 
$\lambda  \in \Lambda$ such that 
$$f(\lambda, \alpha z_{1} + \beta z_{2}) \ \leq \ \alpha 
f(\lambda_{1}, z_{1}) + \beta f(\lambda_{2}, z_{2}) \quad .$$
\label{defimpl}
\end{defi}

In the following  {\bf bipotential convex covers} are defined, as in definition 4.2 \cite{bipo3}.

\begin{defi}   A {\bf bipotential
convex cover} of the non empty set $M$ is a function   
$\displaystyle \lambda \in \Lambda \mapsto b_{\lambda}$ from  $\Lambda$ with 
values in the set  $Bp(X,Y)$, with the 
properties:
\begin{enumerate}
\item[(a)] The set $\Lambda$ is a non empty compact topological space, 
\item[(b)] Let $f: \Lambda \times X \times Y \rightarrow \mathbb{R} \cup
\left\{ + \infty \right\}$ be the function defined by 


$$f(\lambda, x, y) \ = \ b_{\lambda}(x,y) .$$


Then for any $x \in X$ and for any $y \in Y$ the functions 
$f(\cdot, x, \cdot): \Lambda \times Y \rightarrow \bar{\mathbb{R}}$ and 
$f(\cdot, \cdot , y): \Lambda \times X \rightarrow \bar{\mathbb{R}}$ are  lower 
semi continuous  on the product spaces   $\Lambda \times Y$ and respectively 
$\Lambda \times X$ endowed with the standard topology, 
\item[(c)] We have $\displaystyle M  \ = \  \bigcup_{\lambda \in \Lambda} 
M(b_{\lambda})$.
\item[(d)] the functions $f(\cdot, x, \cdot)$ 
and $f(\cdot, \cdot , y)$  are implicitly convex in the sense of Definition 
\ref{defimpl}.
\end{enumerate}
\label{defcover}
\end{defi}

This notion generalizes the one  of a {\bf bi-implicitly 
convex lagrangian cover}.  see Definitions 4.1 and 6.6 \cite{bipo1}. Here we shall give only 
the definition of a convex lagrangian cover, without the bi-implicit convexity hypothesis. 

\begin{defi} Let $M \subset X \times Y$ be a non empty set.  A {\bf convex lagrangian cover }
of  $M$ is a function   $\displaystyle \lambda \in \Lambda \mapsto \phi_{\lambda}$ from  $\Lambda$ with values in the set of lower semicontinuous  and convex functions on $X$, with the 
properties:
\begin{enumerate}
\item[(a)] The set $\Lambda$ is a non empty compact topological space, 
\item[(b)] Let $f: \Lambda \times X \times Y \rightarrow \mathbb{R}$ be the function defined by 
$$f(\lambda, x, y) \ = \ \phi_{\lambda}(x) + \phi^{*}_{\lambda}(y) .$$
Then for any $x \in X$ and for any $y \in Y$ the functions $f(\cdot, x, \cdot)$ and $f(\cdot, \cdot , y)$ are  lower semicontinuous  from  $\Lambda$ with values in the set of lower semicontinuous  and convex functions on $X$, endowed with pointwise convergence topology, 
\item[(c)] we have $$M  \ = \  \bigcup_{\lambda \in \Lambda} M(\phi_{\lambda}) \ .$$
\end{enumerate}
\label{defcover1}
\end{defi}

 A bipotential convex cover  $\displaystyle \lambda \in \Lambda 
\mapsto b_{\lambda}$ such that for any $\lambda \in \Lambda$ the bipotential 
$\displaystyle b_{\lambda}$ is separable is a bi-implicitly 
convex lagrangian cover, see   Definitions 4.1 and 6.6 \cite{bipo1}. . For such covers the sets 
$\displaystyle M(b_{\lambda})$are {\bf maximal cyclically monotone} for any $\lambda \in \Lambda$.

General bipotential convex covers are {\bf not lagrangian}.  (see remark 6.1 
\cite{bipo1} for a justification of the "lagrangian" term). In the language of 
convex analysis this means that in general the sets $\displaystyle M(b_{\lambda})$ are 
not  cyclically monotone. An example is given in section 5 \cite{bipo3}, of a bipotential convex 
cover of the graph of the Coulomb's dry friction law,  which is made by monotone, but not maximally
monotone graphs.

In sections 5 and 8 \cite{bipo1} it is explained that not any BB-graph admits a convex lagrangian cover. 
Moreover, there are   BB-graphs  
admitting (up to reparametrization) only 
one  convex lagrangian cover, as well as BB-graphs 
which have infinitely many lagrangian covers.  The problem of describing 
the set of all convex lagrangian covers of a BB-graph seems to be difficult.

\begin{defi}
Let $\displaystyle \lambda \mapsto b_{\lambda}$ be a bipotential convex cover 
 of the BB-graph $M$. To the cover we associate the function
$b: X \times Y \rightarrow \mathbb{R} \cup \left\{ + \infty \right\}$ 
by the formula 
$$b(x,y) \ = \ \inf \left\{ b_{\lambda}(x,y) 
\mbox{ : } \lambda \in \Lambda\right\} $$
\label{defrecipe}
\end{defi}

We obtained in theorem 4.6 \cite{bipo3} the following result. 

\begin{thm} Let $\displaystyle \lambda \mapsto \phi_{\lambda}$ be a bipotential convex cover  of 
 the BB-graph $M$ and $b: X \times Y \rightarrow R$ the associated function 
 introduced in  Definition
 \ref{defrecipe}. 
Then $b$ is a bipotential and $M=M(b)$. 
\label{thm2}
\end{thm}

In the case of $M=M(\phi)$, with $\phi$ convex and lower semi continuous  
(this corresponds to separable bipotentials), the set $\Lambda$ has only one 
element $\Lambda = \left\{ \lambda \right\}$ 
  and we have only one potential $\displaystyle \phi$. The associated 
bipotential from Definition \ref{defrecipe} is obviously 
$$b(x,y) \ = \ \phi(x) + \phi^{*}(y) \ .$$
This is a bipotential convex cover in a trivial way; the implicit convexity conditions are equivalent with 
the convexity of $\displaystyle \phi$, $\displaystyle \phi^{*}$ respectively.

\section{One more construction result}
\label{newconst}

For simplicity, in this section we shall work only with lower semicontinuous convex functions 
$\phi$ with the property that $\phi \in \Gamma(X,\mathbb{R}) $ and its Fenchel dual 
$\displaystyle \phi^{*} \in \Gamma(Y, \mathbb{R})$. 

We reproduce here 
the following definition of convexity (in a generalized sense), given by 
K. Fan \cite{kyfan} p. 42.

\begin{defi} Let $X$, $Y$ be two  arbitrary non empty sets. The function   
$f:X\times Y  \rightarrow \mathbb{R}$ is {\bf convex} on $X$ {\bf in the sense of Fan} if for any two 
elements $\displaystyle x_{1}, x_{2} \in X$ and for any two numbers 
$\alpha, \beta \in [0,1]$ with $\alpha + \beta = 1$ there exists a $x \in X$ 
such that for all $y \in Y$: 
$$f(x,y) \ \leq \ \alpha f(x_{1},y) + \beta f(x_{2},y) .$$
\label{defkyfan}
\end{defi}

With the help of the previous definition we introduce a new convexity condition 
for a convex lagrangian cover.

\begin{defi}
Let $\displaystyle \lambda \mapsto \phi_{\lambda}$ be a convex lagrangian cover 
of the BB-graph $M$, in the sense of definition \ref{defcover1}. Consider the 
functions: 
$$g: X \times \Lambda  \times X \rightarrow \mathbb{R} \quad , \quad 
h: Y \times \Lambda \times Y \rightarrow \mathbb{R} \quad , $$
given by $\displaystyle g(x,\lambda, z) = \phi_{\lambda}(x) - \phi_{\lambda}(z)$, 
respectively  $\displaystyle h(y, \lambda, u) = \phi_{\lambda}^{*}(y) - 
\phi_{\lambda}^{*}(u)$. 

The cover is {\bf Fan bi-implicitly convex}  if for any $x \in X$, $y \in Y$, 
the functions $g(x, \cdot , \cdot)$, $h(y, \cdot , \cdot)$ are convex in the 
sense of Fan on $\Lambda \times X$, $\Lambda \times Y$ respectively. 
\label{deffic}
\end{defi}

Recall the following  minimax theorem of Fan \cite{kyfan}, Theorem 2. 
In the formulation of the theorem words "convex" and "concave" have the meaning given in definition \ref{defkyfan} 
(more precisely $f$ is concave if $-f$ is convex in the sense of the before 
mentioned definition). 

\begin{thm} (Fan) Let $X$ be a compact Hausdorff space and $Y$ an arbitrary 
set. Let $f$ be a real valued function on $X\times Y$ such that, for every $y 
\in Y$, $f(\cdot ,y)$ is lower semicontinuous on $X$. If $f$ is convex on $X$ 
and concave on $Y$, then we have
$$\min_{x \in X} \sup_{y \in Y} f(x,y) = \sup_{y \in Y} \min_{x \in X} f(x,y) 
\quad . $$
\label{tfan}
\end{thm}

The difficulty of theorem \ref{thm2} boils down to the fact the  class of 
convex functions is not closed with respect to the inf operator. Nevertheless, 
by using Fan theorem \ref{tfan} we get the following general result.

\begin{thm} Let $\Lambda$ be a compact Hausdorff space and 
$\displaystyle \lambda \mapsto \phi_{\lambda} \in \Gamma(X, \mathbb{R})$ be a 
convex lagrangian cover  of the BB-graph $M$, with $\displaystyle
\phi^{*}_{\lambda} \in \Gamma(Y, \mathbb{R})$ for any $\lambda \in \Lambda$,  such that: 
\begin{enumerate}
\item[(a)]   for any $x \in X$ and for any $y \in Y$ the functions 
$\displaystyle \Lambda \ni \lambda  \mapsto \phi_{\lambda}(x) \in \mathbb{R} $ 
and $\displaystyle \Lambda \ni \lambda  \mapsto 
\phi_{\lambda}^{*}(y) \in \mathbb{R}$  are  continuous,  
\item[(b)] the cover is Fan bi-implicitly convex in the sense of definition 
\ref{deffic}.
\end{enumerate}
Then  the function  
$b: X \times Y \rightarrow \mathbb{R}$ defined by
$$b(x,y)  \ = \ \inf \left\{ \phi_{\lambda}(x) + \phi^{*}_{\lambda}(y) \ \mid \  \lambda \in \Lambda \right\}$$
is a bipotential and $M=M(b)$. 
\label{thm4}
\end{thm}

\paragraph{Proof.} 
For some of the details of the proof we refer to the proof of theorem 4.12 
 \ref{thm2}    in  \cite{bipo1}. 
There are five steps in that proof. In order to prove our theorem we have only
 to modify the first two steps: we want to show that for any $x \in  dom(M)$ 
and any $y \in  im(M)$ the functions  $b(\cdot, y)$ and $b(x, \cdot)$ are 
convex and lower semi continuous.
 
For $(x,y) \in X \times Y$ let us define the function 
$\displaystyle \overline{xy}: \Lambda \times X \rightarrow \mathbb{R}$ by 
$$\overline{xy} (\lambda, z) = \langle z,y\rangle + \phi_{\lambda}(x) - \phi_{\lambda}(z)  \quad . $$
We check now that $\displaystyle \overline{xy}$ verifies the hypothesis of 
theorem \ref{tfan}. Indeed, the hypothesis (a) implies that  for any $z \in X$ 
the function  $\displaystyle \overline{xy}(\cdot , z)$  is 
continuous. Notice that 
$$ \overline{xy} (\lambda,z) = \langle z, y \rangle + g(x,\lambda,z) \quad . $$
It follows from hypothesis (b) that the function $\displaystyle \overline{xy}$ is convex on 
$\Lambda$ in the sense of Fan. 

In order to prove the concavity of  
$\displaystyle \overline{xy}$ on $X$, it 
suffices to show that for any $\displaystyle z_{1}, z_{2} \in X$, for any 
$\alpha, \beta \in [0,1]$ such that $\alpha + \beta = 1$, we have the inequality 
$$\overline{xy}(\lambda, \alpha z_{1} + \beta z_{2}) \leq 
\alpha \overline{xy}(\lambda, z_{1}) + \beta \overline{xy}(\lambda, z_{2}) $$ 
for any $\lambda \in \Lambda$.   This inequality  is equivalent with 
$$\langle \alpha z_{1} + \beta z_{2} , y \rangle - \phi_{\lambda}(\alpha z_{1} + 
\beta z_{2}) \leq \alpha \left( \langle z_{1} , y\rangle - \phi_{\lambda}(z_{1} \right) + 
\beta \left( \langle z_{2} , y\rangle - \phi_{\lambda}(z_{2} \right)$$
for any $\lambda \in \Lambda$. But this  is implied by  the 
convexity of $\displaystyle \phi_{\lambda}$ for any $\lambda \in \Lambda$. 

In conclusion the function $\displaystyle \overline{xy}$ satisfies the hypothesis 
of theorem \ref{tfan}.  We deduce that 
\begin{equation}
\min_{\lambda  \in \Lambda} \sup_{z \in X}  \overline{xy}(\lambda, z) = \sup_{z \in X} \min_{\lambda \in \Lambda}\overline{xy}(\lambda, z) \quad . 
\label{eqminsup}
\end{equation}
Let us compute  the two sides of this equality. 

For the left hand side term of (\ref{eqminsup}) we have: 
$$\min_{\lambda  \in \Lambda} \sup_{z \in X}  \overline{xy}(\lambda, z) = \min_{\lambda  \in \Lambda} \sup_{z \in X}  \left\{ \langle z,y\rangle + \phi_{\lambda}(x) - \phi_{\lambda}(z) \right\} = $$
$$= \min_{\lambda  \in \Lambda}  \left\{ \phi_{\lambda}(x) + \sup_{z \in X} \left\{ \langle z,y\rangle - \phi_{\lambda}(z) \right\}  \right\} = $$
$$ =  \min_{\lambda  \in \Lambda}  \left\{ \phi_{\lambda}(x) + \phi^{*}_{\lambda}(y) \right\} = b(x,y) \quad . $$
For the right hand side term of (\ref{eqminsup}) we have: 
$$\sup_{z \in X} \min_{\lambda \in \Lambda}\overline{xy}(\lambda, z) = \sup_{z \in X} \min_{\lambda \in \Lambda}\left\{ \langle z,y\rangle + \phi_{\lambda}(x) - \phi_{\lambda}(z) \right\} = $$
$$= \sup_{z \in X}  \left\{  \langle z,y\rangle -  \max_{\lambda \in \Lambda}\left\{  \phi_{\lambda}(z) - \phi_{\lambda}(x) \right\} \right\} \quad . $$
Let $\displaystyle \overline{x}: X \rightarrow \mathbb{R}$ be the function 
$$ \overline{x}(z) = \max_{\lambda \in \Lambda}\left\{  \phi_{\lambda}(z) - \phi_{\lambda}(x) \right\} \quad . $$
Then the right hand side term of (\ref{eqminsup}) is in fact: 
$$\sup_{z \in X} \min_{\lambda \in \Lambda}\overline{xy}(\lambda, z) = \overline{x}^{*}(y) \quad . $$
Therefore we proved the equality: 
$$b(x,y) =   \overline{x}^{*}(y) \quad . $$
This shows that the function $b$ is convex and lower semicontinuous in the second argument. 

In order to 
prove that $b$ is convex and lower semicontinuous in the first  argument,  
replace $\displaystyle \phi_{\lambda}$ by $\displaystyle 
\phi_{\lambda}^{*}$ in the previous reasoning. \quad $\square$

\vspace{\baselineskip}

\end{document}